# A Proposed Lean Distribution System for Solar Power Plants Using Mathematical Modeling and Simulation Technique

Mohsen Momenitabar
*College of Business, North Dakota State University,*
Fargo, ND 58105, USA
ORCID 0000-0003-2568-1781

Zhila Dehdari Ebrahimi
*College of Business, North Dakota State University,*
Fargo, ND 58105, USA
ORCID 0000-0001-7256-0881

Seyed Hassan Hosseini
*Department of Civil and Industrial Engineering, Sapienza University of Rome*, Lazio 00133, Italy
ORCID 0000-0003-3239-3163

Mohammad Arani
*Department of Systems Engineering, University of Arkansas at Little Rock,*
Little Rock, AR 72204, USA
ORCID 0000-0002-1712-067X

*Abstract*—Today, power waste is one of the most crucial problems which power stations across the world are facing. One of the recent trends of the energy system is the lean management technique. When lean management is indicated by the system, customer value is increased and the waste process in industry or in a power station is reduced. In this paper, first of all, we propose mathematical modeling to reduce the cost of production, and then the simulation technique applies to electricity transmission distribution systems. Furthermore, we consider two criteria for comparison including the different cost of the system and the rate of energy waste during the transmission. The primary approach is to use both of the models in order to draw a comparison between simulation results and mathematical models. Finally, the analysis of the test results done by the CPLEX toolbox of MATLAB Software 2019 leads to a remarkable decrease in the costs of energy demand in the electricity transmitting network distribution system.

## I. Introduction

The development of advanced clean energy technologies worldwide needs to be accelerated urgently in order to diminish the effects of global warming and the threatening climate changes causing fossil fuels. Solar power is thought to be a type of renewable and sustainable energy that is endless [1]–[4]. In practice, solar power is a proven alternative for renewable energies, such as lighting, space and water heating [5]. So solar power stations can be a logical substitute for old plants that lead to a huge amount of atmospheric pollution. Solar plants are primarily used for electricity production that supplies energy for different operations as an essential demand for a modern society. When installing a solar power plant to select the largest electricity supplier with the greatest creation of capability for adding value to the distributed nature of solar electricity generated, it is important to be aware of the photovoltaic capacity of the various cities regions [6]–[8].

Most studies have used the Weibull model, the Rayleigh distribution model, maximum entropy theory, artificial neural networks, geographical information systems, and extract functions algorithms in general for the analysis of wind energy and solar energy [9]. The research also shows how significant solar power plants are the performance of various solar power plants and photovoltaic generations and their financial improvement have been studied in various countries. [10] created a mathematical model, which then analyzed the experimental performance of a solar chemistry system demonstrator. [11] Presented a photovoltaic/thermal (PV/T) design optimization system with both concentrated and concentrated solar radiations. [12] Studied the impact of the generation of high-scale photovoltaic on oscillation of the power system, which specifically considered inter-area oscillation. The effects of photovoltaic on inter-area modes were investigated in New England and New York for different levels of penetration and operational conditions. Increased photovoltaic penetration has been found to have an adverse impact on the critical inter-area mode. The real photovoltaic production from two 100 KWP grid-related systems in the same zone were studied by [13], both of which had the same temperature and radiation fluctuations. Mathematical programming and simulation technique for photovoltaic power and solar power plants were hardly an attempt to be presented as lean production method, though many types of research have been carried out in connection with solar and photovoltaic plants and generations, especially PV. More studies can be fund here [2], [14]–[16].

In this paper, solar and photovoltaic plants and generations of buildings were investigated. First, it is installed in a city that takes into consideration different parameters, including production costs, transition costs, additional power costs, cost of installation in all solar power stages, etc. This is a lean production model consisting of solar panels for each building in the city separately. The models also include the capacity of the solar panels, the number of required panels, the price of all panels, etc. In order to analyze these two model types and to compare them with each other the optimized model is considered to employ the simulation method for both models. Finally, the simulation and mathematical programming results illustrate the advantages and disadvantages of both models presented . The study can also lead users to either use the electricity produced from the construction of solar plants in the city or the photovoltaic generations installed in each building [17]–[19].

## II. SOLAR ENERGY FOR ELECTRICITY PRODUCTION

There are both benefits and drawbacks of using solar energy as a means of power production. Even though there are many benefit factors in using solar energy-like environmental emission reductions, the cost-effectiveness is clean, affordable, endless and easy-to-use compared with expected high oil costs [12]. The initial cost of sun-energy collection equipment is one of the major drawbacks. With the cost of solar panels falling, solar panels are increasingly generating electricity [13]. A solar energy system needs a wide range to host the device in order to provide a source of electricity efficiently. In areas where space is slim or costly, this can be discouraging. Pollution can be detrimental to solar panels as pollution can impair photovoltaic cells' performance.. Clouds may have the same effect because sunlight intensity is partially blocked. The question of the older solar components is even more critical as modern designs include the most severe impact mitigation technologies. In addition, solar power only works when the sun is light, so power needs to be supported. During the night despite the effect of the downside being reduced by the use of solar charge systems, expensive solar systems remain unused, so that the position of solar panels can affect effectiveness. TEP, for example, runs the 5-MW system, one of America's largest solar photovoltaic arrays. But the capacity factor for this generator has been averaged 19 percent for two years , which means that it produces most of the time just 19 percent of its rated output. It is necessary to note that, in order to increase production and reduce cost, solar energy as a renewable energy source is minimal.

The authors in this article have tried using a lean manufacturing approach for mathematical programming and simulation methods on two solar photovoltaic ( PV) and PV generation models to represent the results in the two models.

## III. LEAN MANUFACTURING

Lean manufacturing stresses different metrics like those of waste reduction, performance enhancement and quality improvement. It focuses on cost minimization, optimizing consumer choice, speeding up delivery and improving product and service efficiency. Approximately 8 essential waste forms have been defined by advocates of lean manufacturing, namely excessive manufacturing, waiting for a next phase for the operation, unnecessary material transportation, parts over-processing, stock excess of minimum numbers, unnecessary employee movement and the production and underuse of human capital. Waste can be characterized more generally as any operation that does not add value to the product or service.

The Lean Production Principle is to constantly increase efforts to reduce non-value-adding behavior. Waste typically absorbs energy, but does not produce customers value. The elimination of waste therefore is necessary and can require eliminating the overuse of services or materials. The multidimensional approach focuses on cost reduction, the elimination of value-added activities as well as cellular development. It also includes overall plant management, output smoothing, reduction of installations and related waste disposal [14]. Two cases in this paper are known as modeling. The first is for solar plants built outside of the city, the second is for sleek manufacturing of urban buildings for photovoltaic generations. Consequently, the application of lean output to the second case will exclude energy transfer costs. The cost to produce electricity and land necessary for solar power plants would be subtracted as would transfers and additional processes. The analysis of the two cases shows clearly the optimal model.

## IV. STRUCTURE OF THE MODEL

These formulations are used in this paper for the development of solar and photovoltaic energy plants. For the first and second models, we are trying to reduce and optimize the production costs for solar power plant installed and manufactured in the buildings on the city and on photovoltaic energy. The optimization of the model structure with a combination of lean production and non-linear programming is presented. Finally, the number of solar panels for urban buildings and solar plants used outside the city is displayed, in line with the demand for electricity supply and compared to the third model.

## V. FIRST PROBLEM FORMULATION

In this model, the plant is a photovoltaic power plant with the same output in all solar panels.

### A. Indices

The following indexes are considered:

$j$: index set of power plants ($j = 1,2,…,m$)
$t$: index set of time periods ($t = 1,2,…,T$)

### B. Parameters

The following parameters are considered:

$C_{jt}$: Cost of producing per kilowatt-hour power in power plant j in period $t$

$V_{jt}$: Transferring cost per kilowatt-hour power by power plant j in period $t$

$H_{jt}$: Excess cost of power produced in power plant j in period $t$

$D_j$: Rate of power Consumption from the power plant j (followed by uniform distribution)

$R_j$: The setup cost of power plant j

$\Delta_{jt}^{max}$: Maximum capacity of power plant j in period $t$ in per unit

$\Delta_{jt}^{min}$: Minimum capacity of power plant j in the period $t$ per unit.

$F$: Number of power plants

$NPW_{jt}$: Net present value costs (operational costs) for the power plant j in period $t$

### C. Variables

Variables are defined as followings:

$Y_j$: 1 if the power plant j is used, 0 otherwise.

$Z_{jt}$: Produced power in power plants j per kilowatt-hour in period $t$.

$K_{jt}$: Energy surplus produced by power plants j in period $t$.

$Z_{Optimum}$: Objective Function of the first model.

### D. Model

This model is used to produce optimal variables using a nonlinear programming method.

$$Min\ Z_{optimum} = \sum_{j=1}^{m} R_j Y_j + \sum_{j=1}^{m}\sum_{t=1}^{T}(NPW_{jt}+V_{jt})Y_j Z_j(1+i)^{-t}$$
$$+ \sum_{j=1}^{m}\sum_{t=1}^{T} H_{jt} K_{jt}(1+i)^{-t}$$
(1)

Subject to:
$$\Delta_{jt}^{min} \le Y_j Z_{jt} \le \Delta_{jt}^{max}, \forall j=1,2,...,m, \forall t=1,2,...,T \quad (2)$$
$$K_{jt} = Y_j Z_{jt} - D_{jt}, \forall j=1,2,...,m, \forall t=1,2,...,T \quad (3)$$
$$\sum_{j=1}^{m} Y_j = F, \forall j=1,2,...,m \quad (4)$$
$$Y_j = \{0,1\}, Z_{jt} \ge 0, K_{jt} \in R, \forall j=1,2,...,m,$$
$$\forall t=1,2,...,T \quad (5)$$

In the first part of the objective function, the cost of setting up a Power plant j is considered. In the second part of the objective function, operating costs and transportation costs per kWh of electricity during the period t=1, 2, …, T for the power plant j are shown. Finally, in the third part of the objective function, the cost of additional power generation in the power plant j during period t is indicated. The constraint shown as number two is to prevent exceeding the maximum and minimum amount of electricity production while choosing the power plant. Constraint (3) shows the amount of surplus by subtracting power Consumption from the power plant j. constraint (4) said there are limited number of power plant j. last constraint displays the binary and continuous variables.

## VI. SECOND PROBLEM FORMULATION

This model analyzes photovoltaic power generations based on a mathematical model installed in the buildings in the city. The decision-making processes are presented according to parameters which are essential for the PV analysis.

### A. Parameters

$I$: the interest rate.
$T$: production life (operation period)
$Q$: operation cost
$C$: cost of purchasing one solar panel
$A$: quantity of total consumed criteria
$B$: production capacity of each purchased solar panel
$N$: the number of solar panels will be needed.

### B. Variables

Only variable is defined as followings:
$Z$: electricity produced by each solar panel

### C. Model

In this paragraph, we propose the photovoltaic generation model installed in town buildings. The basis is the mathematical model [20]:

$$F = \frac{A.C.Q.Z}{B} + \int \sum_{t=1}^{T} \frac{A.Q.Z}{B}(1+I)^{-t}.dZ \quad (6)$$

$$F = \frac{A.C.Q.Z}{B} + \int \frac{A.Q.Z}{B}\left(\frac{1-\left(\frac{1}{1+I}\right)^T}{I(I+1)}\right).dZ \quad (7)$$

$$F = \frac{A.C.Q.Z}{B} + \frac{A.Q.Z^2}{2.B}\left(\frac{1-\left(\frac{1}{1+I}\right)^T}{I(I+1)}\right) \quad (8)$$

$$\frac{\partial F}{\partial Z} = 0 \quad (9)$$

$$\frac{A.C.Q}{B} + \frac{A.Q.Z}{B}\left(\frac{1-\left(\frac{1}{1+I}\right)^T}{I(I+1)}\right) = 0 \quad (10)$$

$$\frac{A.Q.Z}{B}\left(\frac{1-\left(\frac{1}{1+I}\right)^T}{I(I+1)}\right) = -\frac{A.C.Q}{B} \quad (11)$$

$$Z^* = \left|\frac{-C.I.(I+1)}{1-\left(\frac{1}{1+I}\right)^T}\right| \quad (12)$$

The model derivative is calculated in order to reach to the optimal amount of $Z^*$. After that, by replacing, the optimized model is proposed:

$$F = \frac{A.C.Q}{B}.\left[\frac{-C.I.(I+1)}{1-\left(\frac{1}{1+I}\right)^T}\right] +$$
$$\frac{A.Q}{2.B}.\left[\frac{-C.I.(I+1)}{1-\left(\frac{1}{1+I}\right)^T}\right]^2.\left(\frac{1-\left(\frac{1}{1+I}\right)^T}{I(I+1)}\right) \quad (13)$$

$$F = \frac{A.C.Q}{B}.\left[\frac{-C.I.(I+1)}{1-\left(\frac{1}{1+I}\right)^T}\right]$$
$$+ \frac{A.Q.C^2}{2.B}.\left[\frac{I(I+1)}{1-\left(\frac{1}{1+I}\right)^T}\right] \quad (14)$$

$$F = \frac{-A.C^2.Q.I(I+1)}{2.B.\left[1-\left(\frac{1}{1+I}\right)^T\right]} \quad (15)$$

## VII. RELATIONSHIP BETWEEN SOLAR PANELS AND POWER PLANTS

In the first model, the following formula calculates $Z_i$:

$$z_i = \frac{-\sum_{j=1}^{m}\int_{a}^{x} H_j.K_j.dK_j(1+i)^{-t} - \sum_{j=1}^{m} R_j.Y_j}{\sum_{j=1}^{m} NPW_j.Y_j} \quad (16)$$

Moreover, Z can be calculated as the first model in the second model:

$$Z = \frac{-N.C.Q + \sqrt{(N.C.Q)^2 + 2.F.N.Q.\left(\frac{1-\left(\frac{1}{1+I}\right)^T}{I(I+1)}\right)}}{N.Q.\left(\frac{1-\left(\frac{1}{1+I}\right)^T}{I(I+1)}\right)} \quad (17)$$

For calculating the number of N*, the two formulas should be equal together as the follow:

$$\frac{-N.C.Q + \sqrt{(N.C.Q)^2 + 2.F.N.Q.\left(\frac{1-\left(\frac{1}{1+I}\right)^T}{I(I+1)}\right)}}{N.Q.\left(\frac{1-\left(\frac{1}{1+I}\right)^T}{I(I+1)}\right)}$$

$$= \frac{-\sum_{j=1}^{m}\int_{a}^{x} H_j.K_j.dK_j(1+i)^{-t} - \sum_{j=1}^{m} R_j.Y_j}{\sum_{j=1}^{m} NPW_j.Y_j} \quad (18)$$

N* is then achieved. It is obvious that the variables that have been continuously obtained by entering the third model parameters so that N* is obtained when Z* is entered.

## VIII. SIMULATION AND THE RESULTS

In this paper, simulation technique is used in order to analyze the solar power plants operating effectively, as it is in some cases extremely complicated and costly to study the physical structures, or even impossible to do so. This method was implemented using the Arena simulation program.

In the first iteration, the primary goal of the simulation technique is to determine the effects of uncertainties on demand. They can be characterized by a certain probability distribution due to their unpredictable demand rates. Using a simulation method, solar power plants are investigated in this segment.

### A. Results of the first model

The results of the optimization model are displayed on three tablets, each with low, medium and high demand, according to different input data. TABLE I. the input data for the optimization model, and the results, on the right, are entered in the MATLAB Software 2019 in toolbox CPLEX program. The best quantity is obtained by $Z = 43,600$ from the objective function.

TABLE I. DATA OF THE FIRST MODEL (LOW DEMAND)

| Row | | $R_j$ | $NPW_{jt}$ | $V_{jt}$ | $H_{jt}$ | $A_{jt}^{min}$ | $A_{jt}^{max}$ | $D_{jt}^{min}$ | $Y_j$ | $Z_{jt}$ | $K_{jt}$ |
|---|---|---|---|---|---|---|---|---|---|---|---|
| j=1 | t=1 | 5*10⁹ | 3.5*10⁴ | 3.4*10⁴ | 1.2*10⁴ | 3*10³ | 4*10³ | 4*10³ | 1 | 4*10³ | 0 |
| j=1 | t=2 | 5*10⁹ | 4.5*10⁴ | 4*10⁴ | 1.3*10⁴ | 3.5*10³ | 4.5*10³ | 4.2*10³ | | 4.2*10³ | 0 |
| j=2 | t=1 | 3.5*10⁹ | 4.5*10⁴ | 3.5*10⁴ | 1.8*10⁴ | 3.2*10³ | 4*10³ | 4.3*10³ | | 4.3*10³ | 0 |
| j=2 | t=2 | 3.5*10⁹ | 2.5*10⁴ | 3*10⁴ | 2*10⁴ | 3.3*10³ | 4*10³ | 4*10³ | 1 | 4*10³ | 0 |
| j=3 | t=1 | 4*10⁹ | 3.2*10⁴ | 2.5*10⁴ | 1.4*10⁴ | 3.1*10³ | 4*10³ | 3.9*10³ | | 0 | 3.9*10³ |
| j=3 | t=2 | 4*10⁹ | 3*10⁴ | 2.9*10⁴ | 1.5*10⁴ | 3.2*10³ | 4*10³ | 3.8*10³ | | 0 | 3.8*10³ |
| j=4 | t=1 | 3*10⁹ | 2.8*10⁴ | 4*10⁴ | 9*10³ | 3.4*10³ | 4*10³ | 4*10³ | 1 | 0 | 0 |
| j=4 | t=2 | 3*10⁹ | 5*10⁴ | 3.9*10⁴ | 1.5*10⁴ | 3.5*10³ | 4.5*10³ | 4.1*10³ | | 4.1*10³ | 0 |

(Input columns on left; Output columns $Z_{jt}$, $K_{jt}$ on right)

TABLE II. Displays the second sample input and output data (medium demand). CPLEX Toolbox calculates the output data and the optimal value for the Z=51.300 objective function. (The first 8 columns are the same as TABLE I. )

TABLE II. INPUT DATA OF THE FIRST MODEL (MEDIUM DEMAND)

| The input | | The output | |
|---|---|---|---|
| $D_{jt}^{min}$ | $Y_j$ | $Z_{jt}$ | $K_{jt}$ |
| 4.5*10³ | 1 | 0 | 5*10³ |
| 4.7*10³ | | 5.2*10³ | 0 |
| 4.8*10³ | 1 | 5.3*10³ | 0 |
| 4.5*10³ | | 0 | 5*10³ |
| 4.4*10³ | 1 | 5.9*10³ | 0 |
| 4.3*10³ | | 5.8*10³ | 0 |
| 4.5*10³ | 0 | 4*10³ | 5*10³ |
| 4.6*10³ | | 4.1*10³ | 5.1*10³ |

In TABLE III. , The third sample data (high demand) is illustrated for input and output. As mentioned earlier, CPLEX Toolbox is available in MATLAB Software with the output data and the optimum amount of the Z=36,000 objective function. (The first 8 columns are the same as TABLE I. )

TABLE III. DATA OF THE FIRST MODEL (HIGH DEMAND)

| Input | Output | | |
|---|---|---|---|
| $D_{jt}^{min}$ | $Y_j$ | $Z_{jt}$ | $K_{jt}$ |
| 5*10³ | 1 | 0 | 4.5*10³ |
| 5.2*10³ | | 4.7*10³ | 4.5*10³ |
| 5.3*10³ | 1 | 4.8*10³ | 0 |
| 5*10³ | | 0 | 0 |
| 4.9*10³ | 1 | 5.4*10³ | 0 |
| 4.8*10³ | | 5.3*10³ | 0 |
| 5*10³ | 0 | 0 | 4.5*10³ |
| 5.1*10³ | | 0 | 4.6*10³ |

B. *Results of the second model*

The results of the second model, based on data, are shown in TABLE IV in this section. that in TABLE IV have five alternatives to calculate photovoltaic energy optimized modeling.

TABLE IV. PARAMETERS OF THE SECOND MODEL.

| | | Alternatives | | | | |
|---|---|---|---|---|---|---|
| | | Korea | China | Taiwan | U.S.A. | Japan |
| Input | I | 0.25 | 0.12 | 0.18 | 0.10 | 0.13 |
| | T (month) | 60 | 12 | 48 | 36 | 24 |
| | Q (year/$) | 80 | 10 | 50 | 90 | 50 |
| | C ($) | 410 | 170 | 250 | 390 | 433 |
| | A (W/year) | 456,250 | 32,850 | 423,400 | 279,225 | 175,200 |
| | B (W) | 250 | 90 | 290 | 255 | 240 |
| Output | Z* | 191.23 | 213.24 | 109.66 | 172.50 | 138.71 |
| | F ($) | 5.7*10⁹ | 6.6*10⁷ | 1*10⁹ | 3.3*10⁹ | 1*10⁹ |

For photovoltaic generation installed on buildings in the city, the optimum number of solar panels used is $1 \times 1029 \times (N^*)$ based on an excellent amount of z obtained from this model.

## IX. CONCLUSION

In this paper, a combination of mathematical modeling and simulation technique is presented for distribution systems in the built and generated solar power generation of electricity transmission system in cities and photovoltaic generation on buildings in town. The contribution of this paper is that the simulation technique and optimization model for solar plants and photovoltaic power generations have not been compared with published data. Both models have been tested and the findings have shown that optimization and simulation models have essentially provided the same findings. Therefore, the two proposed methods described as reliable techniques in order to implement an optimized model with the minimum energy waste and maximum performance, with different solar power plants or photo solar power generations. We have replicated our data and validated the result of previous work [2]. The first and second models are of the same type, as they are very useful in using both methods to analyze the workings of solar plants and the output of photovoltaic energy used in different buildings as reliable techniques. A method used to assess the efficiency of the operational research derived from the two models. Although the results can be obtained in different circumstances from different data from the first and second models, this paper includes a comparison of the solar energy plants Z* and F* and generates better energy generations based on building performance. For example, the solar power plant (low demand) values of $Z$ were $Z = 43,600$ and the photovoltaic power generation value of $F^*$ in Alt1 (Korea) was $F^* = 16,346$. It can also be inferred that solar power plant generations perform better economically than photovoltaic power plants.